\documentclass[12pt]{amsart}

\usepackage{amsmath,amssymb,amsthm,amscd}

\begin{document}

\newcommand{\End}{{\rm{End}\ts}}
\newcommand{\Hom}{{\rm{Hom}}}
\newcommand{\ch}{{\rm{ch}\ts}}
\newcommand{\non}{\nonumber}
\newcommand{\wt}{\widetilde}
\newcommand{\wh}{\widehat}
\newcommand{\ot}{\otimes}
\newcommand{\la}{\lambda}
\newcommand{\La}{\Lambda}
\newcommand{\al}{\alpha}
\newcommand{\be}{\beta}
\newcommand{\ga}{\gamma}
\newcommand{\ka}{\kappa}
\newcommand{\si}{\sigma}
\newcommand{\vp}{\varphi}
\newcommand{\de}{\delta^{}}
\newcommand{\om}{\omega}
\newcommand{\hra}{\hookrightarrow}
\newcommand{\ve}{\varepsilon}
\newcommand{\ts}{\,}
\newcommand{\ee}{e^{}}
\newcommand{\su}{s^{}}
\newcommand{\sii}{\sigma^{-1}}
\newcommand{\vac}{\mathbf{1}}
\newcommand{\di}{\partial}
\newcommand{\qin}{q^{-1}}
\newcommand{\tss}{\hspace{1pt}}
\newcommand{\Sr}{ {\rm S}}
\newcommand{\U}{ {\rm U}}
\newcommand{\Y}{ {\rm Y}}
\newcommand{\AAb}{\mathbb{A}\tss}
\newcommand{\CC}{\mathbb{C}\tss}
\newcommand{\QQ}{\mathbb{Q}\tss}
\newcommand{\SSb}{\mathbb{S}\tss}
\newcommand{\ZZ}{\mathbb{Z}\tss}
\newcommand{\Z}{{\rm Z}}
\newcommand{\Ac}{\mathcal{A}}
\newcommand{\Pc}{\mathcal{P}}
\newcommand{\Qc}{\mathcal{Q}}
\newcommand{\Tc}{\mathcal{T}}
\newcommand{\Sc}{\mathcal{S}}
\newcommand{\Bc}{\mathcal{B}}
\newcommand{\Ec}{\mathcal{E}}
\newcommand{\Hc}{\mathcal{H}}
\newcommand{\Uc}{\mathcal{U}}
\newcommand{\Wc}{\mathcal{W}}
\newcommand{\Ar}{{\rm A}}
\newcommand{\Ir}{{\rm I}}
\newcommand{\Zr}{{\rm Z}}
\newcommand{\gl}{\mathfrak{gl}}
\newcommand{\Pf}{{\rm Pf}}
\newcommand{\oa}{\mathfrak{o}}
\newcommand{\spa}{\mathfrak{sp}}
\newcommand{\g}{\mathfrak{g}}
\newcommand{\z}{\mathfrak{z}}
\newcommand{\Zgot}{\mathfrak{Z}}
\newcommand{\p}{\mathfrak{p}}
\newcommand{\sll}{\mathfrak{sl}}
\newcommand{\agot}{\mathfrak{a}}
\newcommand{\qdet}{ {\rm qdet}\ts}
\newcommand{\cdet}{ {\rm cdet}}
\newcommand{\tr}{ {\rm tr}}
\newcommand{\loc}{{\rm loc}}
\newcommand{\Gr}{ {\rm Gr}\tss}
\newcommand{\sgn}{ {\rm sgn}\ts}
\newcommand{\Sym}{\mathfrak S}
\newcommand{\fand}{\quad\text{and}\quad}
\newcommand{\Fand}{\qquad\text{and}\qquad}
\newcommand{\vt}{{\tss|\hspace{-1.5pt}|\tss}}

\renewcommand{\theequation}{\arabic{section}.\arabic{equation}}

\newtheorem{thm}{Theorem}[section]
\newtheorem{lem}[thm]{Lemma}
\newtheorem{prop}[thm]{Proposition}
\newtheorem{cor}[thm]{Corollary}
\newtheorem{conj}[thm]{Conjecture}

\theoremstyle{definition}
\newtheorem{defin}[thm]{Definition}

\theoremstyle{remark}
\newtheorem{remark}[thm]{Remark}
\newtheorem{example}[thm]{Example}

\newcommand{\bth}{\begin{thm}}
\renewcommand{\eth}{\end{thm}}
\newcommand{\bpr}{\begin{prop}}
\newcommand{\epr}{\end{prop}}
\newcommand{\ble}{\begin{lem}}
\newcommand{\ele}{\end{lem}}
\newcommand{\bco}{\begin{cor}}
\newcommand{\eco}{\end{cor}}
\newcommand{\bde}{\begin{defin}}
\newcommand{\ede}{\end{defin}}
\newcommand{\bex}{\begin{example}}
\newcommand{\eex}{\end{example}}
\newcommand{\bre}{\begin{remark}}
\newcommand{\ere}{\end{remark}}
\newcommand{\bcj}{\begin{conj}}
\newcommand{\ecj}{\end{conj}}

\newcommand{\bal}{\begin{aligned}}
\newcommand{\eal}{\end{aligned}}
\newcommand{\beq}{\begin{equation}}
\newcommand{\eeq}{\end{equation}}
\newcommand{\ben}{\begin{equation*}}
\newcommand{\een}{\end{equation*}}

\newcommand{\bpf}{\begin{proof}}
\newcommand{\epf}{\end{proof}}

\def\beql#1{\begin{equation}\label{#1}}

\title[Fusion procedure]
{Fusion procedure for the Brauer algebra}

\author{A. P. Isaev}
\address{Bogoliubov Laboratory of Theoretical Physics,
Joint Institute for Nuclear Research,
Dubna, Moscow region 141980, Russia}
\email{isaevap@theor.jinr.ru}

\author{A. I. Molev}
\address{School of Mathematics and Statistics,
University of Sydney, NSW 2006, Australia}
\email{alexm\hspace{0.09em}@\hspace{0.1em}maths.usyd.edu.au}

\date{} 


\begin{abstract}
We show that all primitive idempotents
for the Brauer algebra $\Bc_n(\om)$ can be found
by evaluating a rational
function in several variables which has the form
of a product of $R$-matrix type factors.
This provides an analogue of
the fusion procedure for $\Bc_n(\om)$.
\end{abstract}

\maketitle


\section{Introduction}
\label{sec:int}
\setcounter{equation}{0}

It is well known that all primitive idempotents
of the symmetric group $\Sym_n$ can be obtained by taking certain
limit values of the rational function
\beql{phiu}
\Phi(u_1,\dots,u_n)= \prod_{1\leqslant i<j\leqslant n}
\Big(1-\frac{s_{ij}}{u_i-u_j}\Big),
\eeq
where $s_{ij}\in\Sym_n$ is the transposition of $i$ and $j$,
$u_1,\dots,u_n$ are
complex variables and the product is calculated in the group algebra
$\CC[\Sym_n]$ in the lexicographical order on the pairs $(i,j)$.
This construction, which is commonly referred to as the {\it fusion
procedure\/}, goes back to Jucys~\cite{j:yo} and Cherednik~\cite{c:sb}.
Detailed proofs were given by Nazarov~\cite{n:yc}.
A simple version of the fusion procedure
was found in \cite{m:fp}; see also
\cite[Ch.~6]{m:yc} for applications to
the Yangian representation theory and more references.
In more detail, let $T$ be a standard
tableau associated with a partition $\lambda$ of $n$ and let
$c_k=j-i$, if the element $k$ occupies the cell of the tableau in
row $i$ and column $j$. Then the consecutive evaluations
\beql{phieval}
\Phi(u_1,\dots,u_n)\big|_{u_1=c_1}\big|_{u_2=c_2}\dots
\big|_{u_n=c_n}
\eeq
are well-defined and this value yields the
corresponding primitive idempotent $E^{\lambda}_{T}$ multiplied by
the product of the hooks of the diagram of $\lambda$.

In this paper we give a similar
fusion procedure for the Brauer algebra $\Bc_n(\om)$.
This algebra was introduced by Brauer in \cite{b:aw} and its structure
and representation theory was studied by many authors;
see, for instance, Wenzl~\cite{w:sb}, Nazarov~\cite{n:yo},
Leduc and Ram~\cite{lr:rh} and Rui~\cite{r:cs}.
We refer the reader to the review paper by
Barcelo and Ram~\cite{br:cr} for the discussion of
the Brauer algebra in the context of combinatorial
representation theory and more references.
The irreducible representations of $\Bc_n(\om)$
are indexed by all partitions of the nonnegative
integers $n,n-2,n-4,\dots$. If $\la$
is a such partition, then the
{\it updown $\la$-tableaux} $T$
parameterize basis vectors of the corresponding
representation; see Sec.~\ref{sec:bra}.

Consider the rational function
\beql{phiubra}
\Psi(u_1,\dots,u_n)=
\prod_{1\leqslant i<j\leqslant n}
\Big(1-\frac{e_{ij}}{u_i+u_j}\Big)
\prod_{1\leqslant i<j\leqslant n}
\Big(1-\frac{s_{ij}}{u_i-u_j}\Big)
\eeq
with the ordered products as in \eqref{phiu};
the elements $e_{ij},s_{ij}\in\Bc_n(\om)$ are defined in
Sec.~\ref{sec:bra} below.
This function was first introduced by Nazarov~\cite[(3.14)]{n:rt}
in the context of representations
of the classical Lie algebras and twisted Yangians.

Our main result is the following analogue of the fusion
procedure for the Brauer algebra:
given an updown $\la$-tableau $T$,
the consecutive evaluations
\beql{reev}
(u_1-c_1)^{p_1}\dots (u_n-c_n)^{p_n}\ts
\Psi(u_1,\dots,u_n)\big|_{u_1=c_1}\big|_{u_2=c_2}\dots
\big|_{u_n=c_n}
\eeq
are well-defined and this value yields the
corresponding primitive idempotent
$E^{\lambda}_{T}$ multiplied by a nonzero constant
$f(T)$ which is calculated in an explicit form.
Here $p_1,\dots,p_n$ are certain integers depending on $T$
which we call the {\it exponents\/}
of $T$ and the $c_i$ are the {\it contents\/} of $T$; see
Sec.~\ref{sec:bra} for precise definitions.

In the particular case where $\la$ is a partition of $n$,
we thus reproduce some closely
related results of Nazarov~\cite{n:rt};
see, in particular, Propositions~3.2, 3.3 and
formulas~(3.20)--(3.23) there. In fact,
he works with wider
classes of representations of the orthogonal
and symplectic groups $G_N$
parameterized by certain skew Young diagrams
with $n$ boxes. The natural action of $G_N$ in the tensor power
$(\CC^N)^{\ot n}$ commutes with the action
of the Brauer algebra $\Bc_n(\om)$ for a suitably
specialized value of $\om$. Nazarov's formulas
for the idempotents provide remarkable
analogues of the Young symmetrizers
in an explicit form. Their images in $(\CC^N)^{\ot n}$
yield realizations
of the representations of $G_N$ associated with
the skew Young diagrams. Note that the corresponding images
of the factors in \eqref{phiubra} are the values of the
Yang $R$-matrix and its transpose; cf. Remark~\ref{rem:rowcol}
below.

If $\la$ is a partition of $n$, then
all exponents
$p_i$ are equal to zero, while the constant $f(T)$
takes the same value as for \eqref{phieval},
thus making this case quite similar to that of the symmetric group.
The existence of a special monomorphism
$\CC[\Sym_n]\to\Bc_n(\om)$ \cite{bb:sc}
can be regarded as an `explanation' of this analogy.
If $\lambda$ is a partition of $n-2f$ for some $f\geqslant 1$,
then the function \eqref{phiubra} can have zeros or poles
of certain multiplicities at $u_i=c_i$ so that
in place of \eqref{phieval} we need to take
`regularized evaluations' as in \eqref{reev}.

The proof of our main theorem (Theorem~\ref{thm:fus})
follows the approach of \cite{m:fp} and
it is based on the construction of the primitive idempotents
$E^{\lambda}_{T}$ in terms of the
Jucys--Murphy elements for the Brauer algebra.
These elements were introduced independently by
Nazarov~\cite{n:yo} and Leduc and Ram~\cite{lr:rh},
where analogues of Young's seminormal
representations for the Brauer algebra were given.
In a more general context of cellular algebras equipped with
a family of Jucys--Murphy elements
the construction of the primitive idempotents
and seminormal forms was given by Mathas~\cite{m:sf}.

We expect a result similar to Theorem~\ref{thm:fus} to hold
for the Birman--Murakami--Wenzl algebras which will be considered
in our publication elsewhere; cf. \cite{i:qg, imo:ih}.

\section{The Brauer algebra and its representations}
\label{sec:bra}
\setcounter{equation}{0}

Let $n$ be a positive integer and $\om$ an indeterminate.
An $n$-diagram $d$ is a collection
of $2n$ dots arranged into two rows with $n$ dots in each row
connected by $n$ edges such that any dot belongs to only one edge.
The product of two diagrams $d_1$ and $d_2$ is determined by
placing $d_1$ above $d_2$ and identifying the vertices
of the bottom row of $d_1$ with the corresponding
vertices in the top row of $d_2$. Let $s$ be the number of
closed loops obtained in this placement. The product $d_1d_2$ is given by
$\om^{\tss s}$ times the resulting diagram without loops.
The {\it Brauer algebra\/} $\Bc_n(\om)$ is defined as the
$\CC(\om)$-linear span of the $n$-diagrams with the
multiplication defined above.
The dimension of the algebra is $1\cdot 3\cdots (2n-1)$.
The following presentation of $\Bc_n(\om)$ is well-known;
see, e.g., \cite{bw:bl}.

\bpr\label{prop:bradr}
The Brauer algebra $\Bc_n(\om)$ is isomorphic to the algebra
with $2n-2$ generators
$s_1,\dots,s_{n-1},e_1,\dots,e_{n-1}$
and the defining relations
\beql{bradr}
\begin{aligned}
s_i^2&=1,\qquad e_i^2=\om\ts \ee_i,\qquad
\su_i\ee_i=\ee_i\su_i=\ee_i,\qquad i=1,\dots,n-1,\\
s_i s_j &=s_js_i,\qquad \ee_i \ee_j = \ee_j \ee_i,\qquad
\su_i \ee_j = \ee_j\su_i,\qquad
|i-j|>1,\\
s_is_{i+1}s_i&=s_{i+1}s_is_{i+1},\qquad
\ee_i\ee_{i+1}\ee_i=\ee_i,\qquad \ee_{i+1}\ee_i\ee_{i+1}=\ee_{i+1},\\
\su_i\ee_{i+1}\ee_i&=\su_{i+1}\ee_i,\qquad \ee_{i+1}\ee_i
\su_{i+1}=\ee_{i+1}\su_i,\qquad
i=1,\dots,n-2.
\end{aligned}
\non
\end{equation}
\epr

The generators $s_i$ and $\ee_i$ correspond to the following diagrams
respectively:

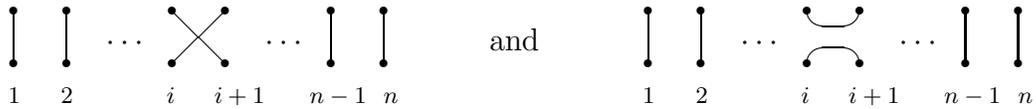
\begin{center}
\begin{picture}(400,60)
\thinlines

\put(10,20){\circle*{3}}
\put(30,20){\circle*{3}}
\put(70,20){\circle*{3}}
\put(90,20){\circle*{3}}
\put(130,20){\circle*{3}}
\put(150,20){\circle*{3}}

\put(10,40){\circle*{3}}
\put(30,40){\circle*{3}}
\put(70,40){\circle*{3}}
\put(90,40){\circle*{3}}
\put(130,40){\circle*{3}}
\put(150,40){\circle*{3}}

\put(10,20){\line(0,1){20}}
\put(30,20){\line(0,1){20}}
\put(70,20){\line(1,1){20}}
\put(90,20){\line(-1,1){20}}
\put(130,20){\line(0,1){20}}
\put(150,20){\line(0,1){20}}

\put(45,25){$\cdots$}
\put(105,25){$\cdots$}

\put(8,5){\scriptsize $1$ }
\put(28,5){\scriptsize $2$ }
\put(68,5){\scriptsize $i$ }
\put(86,5){\scriptsize $i+1$ }
\put(122,5){\scriptsize $n-1$ }
\put(150,5){\scriptsize $n$ }

\put(190,25){\text{and}}

\put(250,20){\circle*{3}}
\put(270,20){\circle*{3}}
\put(310,20){\circle*{3}}
\put(330,20){\circle*{3}}
\put(370,20){\circle*{3}}
\put(390,20){\circle*{3}}

\put(250,40){\circle*{3}}
\put(270,40){\circle*{3}}
\put(310,40){\circle*{3}}
\put(330,40){\circle*{3}}
\put(370,40){\circle*{3}}
\put(390,40){\circle*{3}}

\put(250,20){\line(0,1){20}}
\put(270,20){\line(0,1){20}}
\put(320,20){\oval(20,12)[t]}
\put(320,40){\oval(20,12)[b]}
\put(370,20){\line(0,1){20}}
\put(390,20){\line(0,1){20}}

\put(285,25){$\cdots$}
\put(345,25){$\cdots$}

\put(248,5){\scriptsize $1$ }
\put(268,5){\scriptsize $2$ }
\put(308,5){\scriptsize $i$ }
\put(326,5){\scriptsize $i+1$ }
\put(362,5){\scriptsize $n-1$ }
\put(390,5){\scriptsize $n$ }

\end{picture}
\end{center}

The subalgebra of $\Bc_n(\om)$ generated over $\CC$
by $s_1,\dots,s_{n-1}$
is isomorphic to the group algebra $\CC[\Sym_n]$ so that $s_i$
can be identified with the transposition $(i,i+1)$.
Then for any $1\leqslant i<j\leqslant n$
the transposition $s_{ij}=(i,j)$ can be regarded
as an element of $\Bc_n(\om)$. Moreover, $e_{ij}$
will denote the element of $\Bc_n(\om)$ represented
by the diagram in which the $i$-th and $j$-th dots in the top row,
as well as the $i$-th and $j$-th dots in the bottom
row are connected by an edge, while
the remaining edges
connect the $k$-th dot in the top row
with the $k$-th dot in the bottom row for each $k\ne i,j$.
Equivalently, in terms of the presentation of $\Bc_n(\om)$
provided by Proposition~\ref{prop:bradr},
\ben
s_{ij}=s_i\tss s_{i+1}\dots s_{j-2}\tss
s_{j-1}\tss s_{j-2}\dots s_{i+1}\tss s_i
\Fand
e_{ij}=s_{i,j-1}\tss e_{j-1}\tss s_{i,j-1}.
\een
The Brauer algebra $\Bc_{n-1}(\om)$ can be regarded
as the subalgebra of $\Bc_n(\om)$ spanned by all diagrams
in which the $n$-th dots in the top and bottom rows are
connected by an edge.

The {\it Jucys--Murphy elements\/} $x_1,\dots,x_n$
for the Brauer algebra $\Bc_n(\om)$ were introduced
independently in \cite{lr:rh} and \cite{n:yo};
they are given by the formulas
\ben
x_r=\frac{\om-1}{2}+\sum_{k=1}^{r-1}(s_{kr}-e_{kr}),
\qquad r=1,\dots,n.
\een
The element $x_n$ commutes
with the subalgebra of $\Bc_{n-1}(\om)$. This implies that
the elements $x_1,\dots,x_n$
of $\Bc_n(\om)$ pairwise commute. They can be used
to construct a complete set of pairwise orthogonal primitive
idempotents for the Brauer algebra following
the approach of Jucys~\cite{j:fy} and Murphy~\cite{m:is};
see also \cite{m:sf} for its generalization to a wider class
of cellular algebras. Namely, let $\la$ be a partition of $n-2f$
for some $f\in\{0,1,\dots,\lfloor n/2\rfloor\}$.
We will identify partitions with their diagrams
so that if the parts of $\la$ are $\la_1,\la_2,\dots$ then
the corresponding diagram is
a left-justified array of rows of unit boxes containing
$\la_1$ boxes in the top row, $\la_2$ boxes
in the second row, etc.
The box in row $i$ and column $j$ of a diagram
will be denoted as the pair $(i,j)$.
An {\it updown $\la$-tableau\/} is a sequence
$T=(\La_1,\dots,\La_n)$ of diagrams such that
for each $r=1,\dots,n$ the diagram $\La_r$ is obtained
from $\La_{r-1}$ by adding or removing one box,
where $\La_0=\varnothing$ is the empty diagram and $\La_n=\la$.
To each updown tableau $T$ we attach the corresponding
sequence of {\it contents\/} $(c_1,\dots,c_n)$, $c_r=c_r(T)$,
where
\ben
c_r=\frac{\om-1}{2}+j-i\qquad\text{or}\qquad
c_r=-\Big(\frac{\om-1}{2}+j-i\Big),
\een
if $\La_r$ is obtained by adding the box $(i,j)$ to $\La_{r-1}$
or by removing this box from $\La_{r-1}$,
respectively. The primitive idempotents $E_{T}=E^{\lambda}_{T}$
can now be defined by
the following recurrence formula (we omit the superscripts
indicating the diagrams since they are determined by
the updown tableaux).
Set $\mu=\La_{n-1}$ and consider the updown $\mu$-tableau
$U=(\La_1,\dots,\La_{n-1})$. Let $\alpha$ be the box
which is added to or removed from $\mu$ to get $\la$. Then
\beql{murphyfo}
E_{T}=E_{U}\ts
\frac{(x_n-a_1)\dots (x_n-a_k)}{(c_n-a_1)\dots (c_n-a_k)},
\eeq
where $a_1,\dots,a_k$ are the contents of all boxes
excluding $\alpha$,
which can be removed from or added
to $\mu$ to get a diagram.
When $\la$ runs over all partitions
of $n,n-2,\dots$ and $T$ runs over all updown $\la$-tableaux,
the elements $\{E_{T}\}$ yield a complete set
of pairwise orthogonal primitive
idempotents for $\Bc_n(\om)$. They have the properties
\beql{xiet}
x_r\ts E_{T}=E_{T}\ts x_r=c_r(T)\ts E_{T}, \qquad r=1,\dots,n.
\eeq
Moreover, given an updown tableau $U=(\La_1,\dots,\La_{n-1})$,
we have the relation
\beql{eutet}
E_{U}=\sum_{T} E_{T},
\eeq
summed over all updown tableaux
of the form $T=(\La_1,\dots,\La_{n-1},\La_n)$; we refer
the reader to \cite{lr:rh}, \cite{m:sf} and \cite{n:yo}
for more details.
The relation \eqref{murphyfo} admits the following
equivalent form
\beql{jmform}
E_{T}=E_{U}\ts \frac{u-c_n}{u-x_n}\ts\Big|^{}_{u=c_n},
\eeq
where $u$ is a complex variable. This relation is derived
from \eqref{xiet} and \eqref{eutet} exactly as in the
case of the symmetric group; see \cite{m:fp}.

\section{The fusion procedure}
\label{sec:fus}
\setcounter{equation}{0}

Some combinatorial data extracted from
the updown tableaux will be convenient for
the formulations below. Given an updown $\mu$-tableau
$U=(\La_1,\dots,\La_{n-1})$ we define two infinite matrices
$m(U)$ and $m'(U)$ whose rows and columns are labelled
by positive integers and only a finite number of entries
in each of the matrices is nonzero. The entry $m_{ij}$
of the matrix $m(U)$ (resp., the entry $m'_{ij}$
of the matrix $m'(U)$) equals the number of times
the box $(i,j)$ was added (resp., removed) in the sequence
of diagrams $(\varnothing=\La_0,\La_1,\dots,\La_{n-1})$.
So, the difference $m(U)-m'(U)$ is the matrix
whose all entries are zero except for the $ij$-th matrix
elements equal to $1$ for which
the corresponding boxes $(i,j)$
are contained in the diagram $\mu$.

\bex\label{ex:matmult}
For the updown tableau

\setlength{\unitlength}{0.3em}
\begin{center}
\begin{picture}(80,8)

\put(-10,0){$U=\Big($}

\put(0,0){\line(0,1){2}}
\put(2,0){\line(0,1){2}}
\put(0,0){\line(1,0){2}}
\put(0,2){\line(1,0){2}}
\put(2.5,0){,}

\put(8,0){\line(0,1){2}}
\put(10,0){\line(0,1){2}}
\put(12,0){\line(0,1){2}}
\put(8,0){\line(1,0){4}}
\put(8,2){\line(1,0){4}}
\put(12.5,0){,}

\put(18,-2){\line(0,1){4}}
\put(20,-2){\line(0,1){4}}
\put(22,0){\line(0,1){2}}
\put(18,-2){\line(1,0){2}}
\put(18,0){\line(1,0){4}}
\put(18,2){\line(1,0){4}}
\put(22.5,0){,}

\put(28,-2){\line(0,1){4}}
\put(30,-2){\line(0,1){4}}
\put(28,-2){\line(1,0){2}}
\put(28,0){\line(1,0){2}}
\put(28,2){\line(1,0){2}}
\put(30.5,0){,}

\put(36,0){\line(0,1){2}}
\put(38,0){\line(0,1){2}}
\put(36,0){\line(1,0){2}}
\put(36,2){\line(1,0){2}}
\put(38.5,0){,}

\put(44,-2){\line(0,1){4}}
\put(46,-2){\line(0,1){4}}
\put(44,-2){\line(1,0){2}}
\put(44,0){\line(1,0){2}}
\put(44,2){\line(1,0){2}}
\put(46.5,0){,}

\put(52,-2){\line(0,1){4}}
\put(54,-2){\line(0,1){4}}
\put(56,0){\line(0,1){2}}
\put(52,-2){\line(1,0){2}}
\put(52,0){\line(1,0){4}}
\put(52,2){\line(1,0){4}}
\put(56.5,0){,}

\put(62,-2){\line(0,1){4}}
\put(64,-2){\line(0,1){4}}
\put(66,-2){\line(0,1){4}}
\put(62,-2){\line(1,0){4}}
\put(62,0){\line(1,0){4}}
\put(62,2){\line(1,0){4}}
\put(66.5,0){,}

\put(72,-2){\line(0,1){4}}
\put(74,-2){\line(0,1){4}}
\put(76,0){\line(0,1){2}}
\put(72,-2){\line(1,0){2}}
\put(72,0){\line(1,0){4}}
\put(72,2){\line(1,0){4}}

\put(78,0){$\Big)$}

\end{picture}
\end{center}
\setlength{\unitlength}{1pt}

\bigskip\bigskip
\noindent
the matrices are
\ben
m(U)=\left[\begin{matrix}1&2\\
                     2&1
       \end{matrix}\right]
\Fand
m'(U)=\left[\begin{matrix}0&1\\
                     1&1
       \end{matrix}\right]
\een
where the common zeros in both matrices have been omitted.
\qed
\eex

Furthermore,
for each integer $k$
we define the nonnegative integers $d_k=d_k(U)$ and
$d^{\tss\prime}_k=d^{\tss\prime}_k(U)$
as the respective sums of the entries of the matrices
$m(U)$ and $m'(U)$ on the $k$-th diagonal:
\ben
d_k=\sum_{j-i=k}m_{ij},\qquad d^{\tss\prime}_k=\sum_{j-i=k}m'_{ij}.
\een
So, in Example~\ref{ex:matmult} we have $d_{-1}=d_0=d_1=2$,
while $d^{\tss\prime}_{-1}=d^{\tss\prime}_0=d^{\tss\prime}_1=1$
and the remaining values $d_k$ and $d^{\tss\prime}_k$ are zero.

Finally, for each integer $k$ introduce the parameters
$g_k=g_k(U)$ and $g'_k=g'_k(U)$ by
\beql{gk}
g_k=\de_{k0}+d_{k-1}+d_{k+1}-2\tss d_k,\qquad
g'_k=d^{\tss\prime}_{k-1}+d^{\tss\prime}_{k+1}-2\tss d^{\tss\prime}_k.
\eeq

Now the {\it exponents\/} $p_1,\dots,p_n$ of
an updown $\la$-tableau $T=(\La_1,\dots,\La_n)$ are defined
inductively, so that $p_r$ depends only on the first $r$
diagrams $(\La_1,\dots,\La_r)$ of $T$. Hence, it is sufficient
to define $p_n$. Taking $U=(\La_1,\dots,\La_{n-1})$ we set
\beql{pn}
p_n=1-g_{k_n}(U)\qquad\text{or}\qquad p_n=1-g'_{k_n}(U),
\eeq
respectively, if $\La_n$ is obtained from $\La_{n-1}$
by adding a box on the diagonal $k_n$ or by
removing a box on the diagonal $k_n$.

\bex\label{ex:expon}
The exponents for the updown tableau

\setlength{\unitlength}{0.3em}
\begin{center}
\begin{picture}(50,8)

\put(-10,0){$T=\Big($}

\put(0,0){\line(0,1){2}}
\put(2,0){\line(0,1){2}}
\put(0,0){\line(1,0){2}}
\put(0,2){\line(1,0){2}}
\put(2.5,0){,}

\put(8,0){\line(0,1){2}}
\put(10,0){\line(0,1){2}}
\put(12,0){\line(0,1){2}}
\put(8,0){\line(1,0){4}}
\put(8,2){\line(1,0){4}}
\put(12.5,0){,}

\put(18,-2){\line(0,1){4}}
\put(20,-2){\line(0,1){4}}
\put(22,0){\line(0,1){2}}
\put(18,-2){\line(1,0){2}}
\put(18,0){\line(1,0){4}}
\put(18,2){\line(1,0){4}}
\put(22.5,0){,}

\put(28,-2){\line(0,1){4}}
\put(30,-2){\line(0,1){4}}
\put(28,-2){\line(1,0){2}}
\put(28,0){\line(1,0){2}}
\put(28,2){\line(1,0){2}}
\put(30.5,0){,}

\put(36,0){\line(0,1){2}}
\put(38,0){\line(0,1){2}}
\put(36,0){\line(1,0){2}}
\put(36,2){\line(1,0){2}}
\put(38.5,0){,}

\put(44,-2){\line(0,1){4}}
\put(46,-2){\line(0,1){4}}
\put(44,-2){\line(1,0){2}}
\put(44,0){\line(1,0){2}}
\put(44,2){\line(1,0){2}}
\put(46.5,0)

\put(48,0){$\Big)$}

\end{picture}
\end{center}
\setlength{\unitlength}{1pt}

\bigskip\bigskip
\noindent
are $p_1=p_2=p_3=0$, $p_4=p_5=1$, $p_6=2$.
\qed
\eex

The constants $f(T)$ which we mentioned in the Introduction
are defined inductively by the formula
\beql{ft}
f(T)=f(U)\ts\vp(U,T),
\eeq
where $U=(\La_1,\dots,\La_{n-1})$ and
$T=(\La_1,\dots,\La_n)$. Here
\ben
\vp(U,T)=\prod_{k\ne k_n}(k_n-k)^{g_k}
\prod_k(k_n+k+\om-1)^{g^{\tss\prime}_k}
\een
or
\ben
\vp(U,T)=\prod_{k\ne k_n}(-k_n+k)^{g^{\tss\prime}_k}
\prod_k(-k_n-k-\om+1)^{g_k},
\een
if $\La_n$ is obtained from $\La_{n-1}$
by adding or removing a box on the diagonal $k_n$,
respectively, where the products are taken
over all integers $k$, while $g_k=g_k(U)$ and $g'_k=g'_k(U)$.
Note that only a finite number of the parameters
$g_k$ and $g'_k$ are nonzero so that
each product in the above formulas contains only a finite number
of factors not equal to $1$.

\bpr\label{prop:symmgr}
If $T=(\La_1,\dots,\La_n)$
is an updown $\la$-tableau and $\la$ is a partition of $n$, then
all exponents $p_1,\dots,p_n$ of $T$ are equal to zero, while
$f(T)$ equals the product of the hooks of $\la$.
\epr

\bpf
Set $U=(\La_1,\dots,\La_{n-1})$ and $\mu=\La_{n-1}$. The nonzero
entries of the matrix $m(U)$ are
equal to $1$; these are the $ij$-th matrix elements such that
the corresponding boxes $(i,j)$ are contained
in the diagram $\mu$. Furthermore,
all entries of the matrix $m'(U)$ are zero. Hence,
the parameters $g'_k(U)$ are all zero, while the nonzero values
of $g_k(U)$ are equal to $\pm1$. The value $1$ (resp., $-1$)
corresponds to
those diagonals $k$ where a box can be added to
(resp., removed from) the diagram $\mu$.
This proves that $p_r=0$ for all $r$ and
the claim about $f(T)$ is also easily verified.
\epf

Consider now the rational function $\Psi(u_1,\dots,u_n)$
with values in the Brauer algebra $\Bc_n(\om)$
defined by \eqref{phiubra}.
We can now prove our main theorem.

\bth\label{thm:fus}
For any updown tableau $T=(\La_1,\dots,\La_n)$
the consecutive evaluations
\ben
(u_1-c_1)^{p_1}\dots (u_n-c_n)^{p_n}\ts
\Psi(u_1,\dots,u_n)\big|_{u_1=c_1}\big|_{u_2=c_2}\dots
\big|_{u_n=c_n}
\een
are well-defined. The corresponding value coincides with
$f(T)\tss E_{T}$.
\eth

\bpf
The proof of the theorem
will follow from a sequence of lemmas.

\ble\label{lem:ybe}
The function $\Psi(u_1,\dots,u_n)$ can be written
in the equivalent form
\begin{multline}\label{psedec}
\Psi(u_1,\dots,u_n)\\
=\prod_{r=2,\dots,n}^{\longrightarrow}
\Big(1-\frac{e_{r-1,r}}{u_{r-1}+u_r}\Big)\dots
\Big(1-\frac{e_{1,r}}{u_{1}+u_r}\Big)
\Big(1-\frac{s_{1,r}}{u_1-u_r}\Big)\dots
\Big(1-\frac{s_{r-1,r}}{u_{r-1}-u_r}\Big),
\end{multline}
where the factors are ordered in accordance with
the increasing values of $r$.
\ele

\bpf
This follows by using the easily verified identities
for the rational functions in $u$ and $v$
with values in $\Bc_n(\om)$:
if $i<j<r$ then
\beql{ybetr}
\Big(1-\frac{e_{ir}}{u}\Big)
\Big(1-\frac{e_{jr}}{v}\Big)
\Big(1-\frac{s_{ij}}{u-v}\Big)=
\Big(1-\frac{s_{ij}}{u-v}\Big)
\Big(1-\frac{e_{jr}}{v}\Big)
\Big(1-\frac{e_{ir}}{u}\Big).
\eeq
If the indices $i,j,k,l$ are distinct, then
the elements $e_{ij}$ and $e_{kl}$ of $\Bc_n(\om)$ commute.
Therefore, we can represent the first product occurring
in \eqref{phiubra} as
\begin{multline}\non
\prod_{1\leqslant i<j\leqslant n}
\Big(1-\frac{e_{ij}}{u_i+u_j}\Big)=
\prod_{1\leqslant i<j\leqslant n-1}
\Big(1-\frac{e_{ij}}{u_i+u_j}\Big)\\
{}\times{}\Big(1-\frac{e_{1,n}}{u_{1}+u_n}\Big)\dots
\Big(1-\frac{e_{n-1,n}}{u_{n-1}+u_n}\Big).
\end{multline}
Now, using the identities \eqref{ybetr} repeatedly,
we get
\begin{multline}\non
\Big(1-\frac{e_{1,n}}{u_{1}+u_n}\Big)\dots
\Big(1-\frac{e_{n-1,n}}{u_{n-1}+u_n}\Big)
\prod_{1\leqslant i<j\leqslant n-1}
\Big(1-\frac{s_{ij}}{u_i-u_j}\Big)\\
{}=\prod_{1\leqslant i<j\leqslant n-1}
\Big(1-\frac{s_{ij}}{u_i-u_j}\Big)
\Big(1-\frac{e_{n-1,n}}{u_{n-1}+u_n}\Big)
\dots\Big(1-\frac{e_{1,n}}{u_{1}+u_n}\Big).
\end{multline}
Hence the function \eqref{phiubra} can be written as
\begin{multline}\label{indude}
\Psi(u_1,\dots,u_n)
=\Psi(u_1,\dots,u_{n-1})\\
{}\times{}\Big(1-\frac{e_{n-1,n}}{u_{n-1}+u_n}\Big)
\dots\Big(1-\frac{e_{1,n}}{u_{1}+u_n}\Big)
\Big(1-\frac{s_{1,n}}{u_1-u_n}\Big)
\dots\Big(1-\frac{s_{n-1,n}}{u_{n-1}-u_n}\Big),
\end{multline}
and the decomposition \eqref{psedec} follows by
the induction on $n$.
\epf

Lemma~\ref{lem:ybe} allows us to use the induction on
$n$ to prove the theorem.
By the induction hypothesis, setting $u=u_n$ we get
\begin{multline}\label{induge}
(u_1-c_1)^{p_1}\dots (u_n-c_n)^{p_n}\ts
\Psi(u_1,\dots,u_n)\big|_{u_1=c_1}\big|_{u_2=c_2}\dots
\big|_{u_{n-1}=c_{n-1}}\\
{}=f(U)\tss E_U\ts (u-c_n)^{p_n}
\Big(1-\frac{e_{n-1,n}}{c_{n-1}+u}\Big)
\dots\Big(1-\frac{e_{1,n}}{c_{1}+u}\Big)
\Big(1-\frac{s_{1,n}}{c_1-u}\Big)
\dots\Big(1-\frac{s_{n-1,n}}{c_{n-1}-u}\Big),
\end{multline}
where $U$ is the updown tableau $(\La_1,\dots,\La_{n-1})$.
The next lemma will allow us to simplify this expression.

\ble\label{lem:jmsi}
We have the identity
\begin{multline}\label{jm}
E_U\ts
\Big(1-\frac{e_{n-1,n}}{c_{n-1}+u}\Big)
\dots\Big(1-\frac{e_{1,n}}{c_{1}+u}\Big)
\Big(1-\frac{s_{1,n}}{c_1-u}\Big)
\dots\Big(1-\frac{s_{n-1,n}}{c_{n-1}-u}\Big)\\
{}=\frac{u-c_1}{u-c_n}\ts
\prod_{r=1}^{n-1}\Big(1-\frac{1}{(u-c_r)^2}\Big)\ts
E_U\ts \frac{u-c_n}{u-x_n}.
\end{multline}
\ele

\bpf
Note that the Jucys--Murphy element $x_n$
commutes with $E_U$, and the inverses of the expressions
occurring in the product are found by
\ben
\Big(1-\frac{s_{r,n}}{c_{r}-u}\Big)^{-1}
\Big(1-\frac{1}{(u-c_r)^2}\Big)=
\Big(1+\frac{s_{r,n}}{c_{r}-u}\Big)
\een
and
\ben
\Big(1-\frac{e_{r,n}}{c_{r}+u}\Big)^{-1}
=\Big(1+\frac{e_{r,n}}{c_{r}+u-\om}\Big),
\een
where we have used the relations $s_{r,n}^2=1$ and
$e_{r,n}^2=\om\ts e_{r,n}$. Hence,
relation \eqref{jm} is equivalent to
\beql{indsi}
\bal
&E_U\ts \Big(1+\frac{s_{n-1,n}}{c_{n-1}-u}\Big)
\dots \Big(1+\frac{s_{1,n}}{c_{1}-u}\Big)
\Big(1+\frac{e_{1,n}}{c_{1}+u-\om}\Big)\dots
\Big(1+\frac{e_{n-1,n}}{c_{n-1}+u-\om}\Big)\\
{}={}&E_U\ts \frac{u-x_n}{u-c_1}.
\eal
\eeq
We embed the Brauer algebra $\Bc_n(\om)$ into $\Bc_m(\om)$
for some $m\geqslant n$ and
verify by induction on $n$
a more general identity
\beql{indsige}
\bal
&E_U\ts \Big(1+\frac{s_{n-1,m}}{c_{n-1}-u}\Big)
\dots \Big(1+\frac{s_{1,m}}{c_{1}-u}\Big)
\Big(1+\frac{e_{1,m}}{c_{1}+u-\om}\Big)\dots
\Big(1+\frac{e_{n-1,m}}{c_{n-1}+u-\om}\Big)\\
{}={}&E_U\ts \frac{u-x^{(m)}_n}{u-c_1},
\eal
\eeq
where
\ben
x^{(m)}_n=\frac{\om-1}{2}+\sum_{k=1}^{n-1}(s_{km}-e_{km}).
\een
By \eqref{eutet} we
have $E_U=E_U\ts E_W$,
where $W$ is the updown tableau $(\La_1,\dots,\La_{n-2})$.
Hence, using the induction hypothesis
we can write the left hand side of \eqref{indsige} as
\ben
\bal
&E_U\ts \Big(1+\frac{s_{n-1,m}}{c_{n-1}-u}\Big)\ts
E_W\ts \frac{u-x^{(m)}_{n-1}}{u-c_1}\ts
\Big(1+\frac{e_{n-1,m}}{c_{n-1}+u-\om}\Big)
=\frac{1}{u-c_1}\ts
{}E_{U}\\
{}&\times\Big(u-x_{n-1}^{(m)} +
\frac{s_{n-1,m}\tss(u-x_{n-1}^{(m)})}{c_{n-1}-u}
+\frac{(u-x_{n-1}^{(m)})\tss e_{n-1,m}}{c_{n-1}+u-\om}
+\frac{s_{n-1,m}\tss(u-x_{n-1}^{(m)})\tss e_{n-1,m}}
{(c_{n-1}-u)(c_{n-1} +u-\om)}\Big).
\eal
\een
Now we use the following relations in $\Bc_m(\om)$
which hold for
$1\leqslant r<n-1$:
\ben
s_{n-1,m}\tss s_{r,m}=s_{r,n-1}\tss s_{n-1,m},\qquad
s_{n-1,m}\tss e_{r,m}=e_{r,n-1}\tss s_{n-1,m}
\een
and
\ben
s_{r,m}\tss e_{n-1,m}=e_{r,n-1}\tss e_{n-1,m},\qquad
e_{r,m}\tss e_{n-1,m}=s_{r,n-1}\tss e_{n-1,m}.
\een
They imply that
\ben
s_{n-1,m}\tss x_{n-1}^{(m)}=x_{n-1}\tss s_{n-1,m}
\een
and
\ben
x_{n-1}^{(m)}\tss e_{n-1,m}=\big(\om-1-x_{n-1}\big)\tss e_{n-1,m}.
\een
Together with the relation
$E_U\ts x_{n-1}=c_{n-1}\ts E_U$ implied by \eqref{xiet},
this allows us to bring the left hand side of \eqref{indsige}
to the form
\ben
\frac{1}{u-c_1}\ts E_{U}
\Big(u-x_{n-1}^{(m)}-s_{n-1,m}+e_{n-1,m}\Big)=
E_{U}\ts\frac{u-x_{n}^{(m)}}{u-c_1},
\een
as required.
\epf

Due to Lemma~\ref{lem:jmsi}, in order to complete
the proof of the theorem,
we need to show that the rational function
\ben
f(U)\tss
(u-c_1)\ts
\prod_{r=1}^{n-1}\Big(1-\frac{1}{(u-c_r)^2}\Big)
\ts(u-c_n)^{p_n-1}\cdot
E_U\ts \frac{u-c_n}{u-x_n}
\een
is regular at $u=c_n$ and its value equals $f(T)\tss E_T$.
Using the parameters \eqref{gk}, we can write
this expression as
\ben
f(U)\ts\prod_{k}\Big(u-\frac{\om-1}{2}-k\Big)^{g_k}
\prod_{k}\Big(u+\frac{\om-1}{2}+k\Big)^{g^{\tss\prime}_k}\ts
(u-c_n)^{p_n-1}\cdot
E_U\ts \frac{u-c_n}{u-x_n},
\een
where $k$ runs over the set of integers.
If the diagram $\La_n$ is obtained from $\La_{n-1}$
by adding or removing a box on the diagonal $k_n$, then
the value of the content $c_n$ is given by the respective formulas
\ben
c_n=\frac{\om-1}{2}+k_n\qquad\text{or}\qquad
c_n=-\Big(\frac{\om-1}{2}+k_n\Big).
\een
The definition of the exponents \eqref{pn},
and the constants $f(T)$ in \eqref{ft} together with
\eqref{jmform} imply the desired statement.
\epf

The following corollary is immediate from
Proposition~\ref{prop:symmgr} and Theorem~\ref{thm:fus};
cf. \cite{m:fp}, \cite{n:rt}.

\bco\label{cor:sym}
If $T=(\La_1,\dots,\La_n)$
is an updown $\la$-tableau and $\la$ is a partition of $n$,
then the consecutive evaluations
\ben
\Psi(u_1,\dots,u_n)\big|_{u_1=c_1}\big|_{u_2=c_2}\dots
\big|_{u_n=c_n}
\een
are well-defined. The corresponding value coincides with
$H(\la)\tss E_{T}$, where $H(\la)$ is the product of
the hooks of $\la$.
\qed
\eco

\bre\label{rem:rowcol}
In two particular cases where $\la$ is a row- or column-diagram
with $n$ boxes, one can write alternative multiplicative
expressions associated with the respective tableaux.
Namely, the primitive idempotent corresponding
to the only updown $(n)$-tableau is proportional to
\ben
\prod_{1\leqslant i<j\leqslant n}
\Big(1+\frac{s_{ij}}{j-i}-\frac{e_{ij}}{j-i+\om/2-1}\Big),
\een
while the primitive idempotent corresponding
to the updown $(1^n)$-tableau is proportional to
\ben
\prod_{1\leqslant i<j\leqslant n}
\Big(1-\frac{s_{ij}}{j-i}\Big),
\een
with both products taken in the lexicographical order
on the pairs $(i,j)$. These formulas are easily
verified by using the well-known fact that
the rational function
\ben
R_{ij}(u)=1-\frac{s_{ij}}{u}+\frac{e_{ij}}{u-\om/2+1}
\een
is a solution of the Yang--Baxter equation
\ben
R_{12}(u)\ts R_{13}(u+v)\ts R_{23}(v)=
R_{23}(v)\ts R_{13}(u+v)\ts R_{12}(u);
\een
see \cite{zz:rf}.
These multiplicative formulas for the idempotents
do not seem to have natural
analogues for general updown tableaux. Note, however, that
the following alternative
rational function in the case of $\Bc_3(\om)$
can be used instead of $\Psi(u_1,u_2,u_3)$
in the formulation of the fusion procedure:
\begin{multline}\non
\wt\Psi(u_1,u_2,u_3) = \Big(1 - (u_1-u_2)\tss s_1
+ \frac{u_1-u_2-1}{u_1+u_2} \ts e_1   \Big)\\
{}\times{}\Big(1 -  (u_1-u_3)\tss s_2 +
\frac{u_1-u_3-2}{u_2+u_3} \ts e_{2}
\Big)
\Big(1 - (u_1-u_2)\tss s_1  + \frac{u_1-u_2-1}{u_1+u_2}
\ts e_1   \Big).
\end{multline}
\ere

\section*{Acknowledgments}

We are grateful to Maxim Nazarov and Oleg Ogievetsky for
valuable discussions.
We acknowledge the support of the Australian Research
Council. The work of the first author was supported by the grants
RFBR 08-01-00392-a, RFBR-CNRS 07-02-92166-a
and RF Grant N.Sh. 195.2008.2. He would like to thank
the School of Mathematics and Statistics
of the University of Sydney for the warm
hospitality during his visit.

\end{document}